\documentclass[twoside,onecolumn,centertags,final,a4paper,12pt,reqno,fleqn]{amsart}
\usepackage[left=1.5cm, bottom=3cm]{geometry}
\usepackage{cite}
\usepackage{amsmath,amssymb,amsfonts}
\usepackage{amsthm}
\usepackage{algorithm,algcompatible,amsmath}
\algnewcommand\INPUT{\item[\textbf{Input:}]}%
\algnewcommand\OUTPUT{\item[\textbf{Output:}]}%
\usepackage{graphicx}
\usepackage{textcomp}
\usepackage{caption}
\usepackage{subcaption}
\usepackage{array}
\usepackage{multirow}
\usepackage{adjustbox}
\newcolumntype{C}[1]{>{\centering\arraybackslash}p{#1}}

\usepackage{enumitem}

\newtheorem{theorem}{Theorem}

\newtheorem{example}[theorem]{Example}

\theoremstyle{remark}

\newcolumntype{C}[1]{>{\centering\arraybackslash}m{#1}}
\usepackage{booktabs}
\usepackage{float}
\usepackage{hyperref}
\usepackage{xcolor}

\allowdisplaybreaks

\begin{document}
	\title[]{Tropical Geometry Based Edge Detection Using Min-Plus and Max-Plus Algebra}
		\author[]{SHIVAM KUMAR JHA S$^{1*}$ and Jaya NN Iyer$^{1}$ \\
	 $^{1}$Department of Mathematics,\\
	  The Institute of Mathematical Sciences (HBNI),\\
		CIT Campus, Taramani, Chennai - 600113, Tamil Nadu, India\\
	 E\lowercase{mail}:\lowercase{$^{1*}$spshivamjha@gmail.com, $^{1}$jniyer@imsc.res.in}\\
 $^{*}$C\lowercase{orresponding Author}}

	
	\begin{abstract}
This paper proposes a tropical geometry-based edge detection framework that reformulates convolution and gradient computations using min-plus and max-plus algebra. The tropical formulation emphasizes dominant intensity variations, contributing to sharper and more continuous edge representations. Three variants are explored: an adaptive threshold-based method, a multi-kernel min-plus method, and a max-plus method emphasizing structural continuity. The framework integrates multi-scale processing, Hessian filtering, and wavelet shrinkage to enhance edge transitions while maintaining computational efficiency. Experiments on MATLAB built-in grayscale and color images suggest that tropical formulations integrated with classical operators, such as Canny and LoG, can improve boundary detection in low-contrast and textured regions. Quantitative evaluation using standard edge metrics indicates favorable edge clarity and structural coherence. These results highlight the potential of tropical algebra as a scalable and noise-aware formulation for edge detection in practical image analysis tasks.
 	\end{abstract}
	
	\maketitle

\section{Introduction}

Edge detection is a fundamental technique in image processing, crucial for tasks such as object recognition, image segmentation, and feature extraction \cite{Ganesan2010}. It involves identifying intensity or color transitions that define object boundaries within an image. Classical methods, including the Canny edge detector, Laplacian of Gaussian (LoG) \cite{Mallat1989}, and gradient-based operators like Sobel \cite{Sobel1968,Duda1973}, Prewitt \cite{Prewitt1970}, and Roberts \cite{Marr1980}, each offer trade-offs in terms of noise sensitivity, structural accuracy, and computational complexity.

Tropical Geometry (TG), a combinatorial framework derived from algebraic geometry \cite{CuninghameGreen1979}, presents an alternative formulation for edge detection. Rather than relying on traditional arithmetic or differential operators, TG employs min-plus and max-plus algebra to process image intensity variations. This approach captures dominant transitions while suppressing minor fluctuations, offering potential advantages in robustness and structure preservation. Related work, such as the Homotopy Parametric Value (LHP) method \cite{Jha24}, demonstrates how algebraic structures can support piecewise linear approximations and valuation-based transformations in imaging.

The proposed framework utilizes tropical algebra to redefine convolution and gradient computations. Image intensity values are mapped into a tropical polynomial space, where edge information is extracted through convex structures and tropical intersections. This formulation is applied across multiple directional kernels and integrated with multi-scale processing, wavelet shrinkage, and Hessian filtering to improve edge continuity. Comparative evaluations with conventional methods—including Canny, Sobel, Prewitt, Roberts, and LoG—are conducted on MATLAB built-in images, both grayscale and color.

A robust edge detection strategy also depends on statistical properties of the image. Key quality factors such as contrast \cite{Johnson1990}, correlation, energy, entropy \cite{Susanj2015}, and homogeneity play an important role in distinguishing meaningful edges from background noise or texture. These measures are used not only in preprocessing and enhancement but also in evaluating algorithmic performance \cite{Li2024}. For example, contrast-based methods enhance visibility in low-light conditions, while entropy captures spatial uncertainty that may obscure fine boundaries. Homogeneity contributes by preserving consistent regions while attenuating spurious gradients. A broader survey of edge detection strategies \cite{Sun2022} provides context for assessing such quality metrics in contemporary algorithms.

This study focuses on three tropical variants designed to demonstrate the flexibility and adaptability of the tropical framework. The $TG^{\text{adapt tresh}}_{min-plus}$ method combines tropical gradients with adaptive thresholding for high-contrast edge delineation. The $TG^{\text{8 kernels}}_{min-plus}$ variant applies multi-directional kernel fusion to capture fine structural details, while $TG^{\text{4 kernels}}_{max-plus}$ emphasizes edge continuity and spatial coherence through max-based tropical convolution. These variants are evaluated in terms of both qualitative visual output and quantitative metrics, including Edge Mean Entropy (EME), contrast, and entropy.

The tropical formulation is also applied as a structural extension to classical operators. For instance, the Canny detector \cite{Canny1986}—a widely used method that integrates smoothing, gradient estimation, non-maximum suppression, and hysteresis thresholding—can be adapted to tropical settings. Here, Gaussian smoothing is replaced with tropical polynomial filtering, while gradients are computed through tropical min-plus algebra. Edge thinning and continuity are guided by tropical convex hull representations and adaptive hysteresis in the tropical domain.

This pipeline involves the following stages:
\begin{enumerate}
    \item Tropical Smoothing: Gaussian blur is replaced by a structure-preserving tropical polynomial transformation.
    \item Tropical Gradient Computation: Gradients are derived via min-plus algebra for improved noise tolerance.
    \item Tropical Non-Maximum Suppression: Convex geometry is used to refine edge maps.
    \item Adaptive Hysteresis Thresholding: A two-level threshold strategy adapted to the tropical domain enhances edge connectivity.
\end{enumerate}

Previous studies modifying the Canny detector \cite{Zhang2021} have proposed improvements in specific domains. The tropical formulation aims to extend this by improving edge continuity and resilience to noise, particularly in complex or low-contrast scenes.

Traditional gradient-based methods like Roberts \cite{Marr1980}, Prewitt \cite{Prewitt1970}, and Sobel \cite{Sobel1968} utilize small convolution kernels to approximate image gradients in horizontal (\(P_x\)) and vertical (\(P_y\)) directions. For example:

\[
R_x =
\begin{bmatrix}
1 & 0 \\
0 & -1
\end{bmatrix},
\quad
R_y =
\begin{bmatrix}
0 & 1 \\
-1 & 0
\end{bmatrix}
\]

\[
P_x =
\begin{bmatrix}
1 & 1 & 1 \\
0 & 0 & 0 \\
-1 & -1 & -1
\end{bmatrix},
\quad
P_y =
\begin{bmatrix}
-1 & 0 & 1 \\
-1 & 0 & 1 \\
-1 & 0 & 1
\end{bmatrix}
\]

\[
S_x =
\begin{bmatrix}
-1 & -2 & -1 \\
0 & 0 & 0 \\
1 & 2 & 1
\end{bmatrix},
\quad
S_y = S_x^T.
\]

Their corresponding gradient magnitudes are given by:

\[
\mu_R = \sqrt{R_x^2 + R_y^2}, \quad
\mu_P = \sqrt{P_x^2 + P_y^2}, \quad
\mu_S = \sqrt{S_x^2 + S_y^2}.
\]

While effective for edge approximation, these operators are sensitive to noise and illumination changes. The proposed tropical gradient-based methods address this by computing gradients through tropical algebra, avoiding additive noise accumulation.

\begin{itemize}
    \item Tropical Convolution: Gradient estimation is performed using min-plus and max-plus operations.
    \item Tropical Edge Enhancement: Edge strength is amplified through tropical maxima, reducing noise influence.
    \item Tropical LoG Integration: The LoG method is reformulated with tropical polynomials to emphasize significant contrasts and suppress background noise.
\end{itemize}

These adaptations enhance edge sharpness, local contrast, and resilience to distortions. Experimental results and comparative analyses support the potential of tropical gradient-based methods as a complementary alternative to classical techniques.

\section{Methodology}
The proposed framework implements edge detection by formulating convolution and gradient operations within the algebraic structure of tropical geometry. It employs min-plus and max-plus algebra to construct directional kernels that emphasize dominant intensity transitions while suppressing noise. The framework integrates three key components: tropical convolution using multi-directional kernels, adaptive thresholding for edge refinement, and wavelet-based enhancement for preserving fine details. These components are applied across both grayscale and color images, producing edge maps that maintain structural continuity and are resilient to contrast variations. The methodology is designed to be modular, allowing integration with conventional operators such as Canny and LoG while preserving the core advantages of tropical algebra. 

Convolution, a core operation in image processing, is defined as  

\begin{equation}  
O(x,y) = \sum_{i,j} K(i,j) \cdot I(x-i, y-j),  
\end{equation}  

where \( O(x,y) \) is the output pixel, \( K(i,j) \) represents the convolution kernel, and \( I(x-i, y-j) \) represents the intensity of the input image at neighboring positions. This operation aggregates intensity values using summation and multiplication. In the tropical setting, convolution is reformulated by replacing summation with the min (or max) operation and multiplication with addition, leading to  

\begin{equation}  
O_T(x,y) = \min_{i,j} \left[ K(i,j) + I(x-i, y-j) \right].  
\end{equation}  

\begin{table}
   \centering
   \resizebox{\textwidth}{!}{%
  \begin{tabular}{cccc}
    
    \subcaptionbox{0$^{\circ}$}{$\begin{bmatrix}
      1 & 1 & 1 \\
      0 & 0 & 0 \\
     -1 & -1 & -1 
    \end{bmatrix}$} & 
    
    \subcaptionbox{90$^{\circ}$}{$\begin{bmatrix}
      1 & 0 & -1 \\
      1 & 0 & -1 \\
      1 & 0 & -1 
    \end{bmatrix}$} & 
    
    \subcaptionbox{180$^{\circ}$}{$\begin{bmatrix}
     -1 & -1 & -1 \\
      0 &  0 &  0 \\
      1 &  1 &  1 
    \end{bmatrix}$} & 
    
    \subcaptionbox{270$^{\circ}$}{$\begin{bmatrix}
     -1 & 0 & 1 \\
     -1 & 0 & 1 \\
     -1 & 0 & 1 
    \end{bmatrix}$} \\

    \subcaptionbox{45$^{\circ}$}{$\begin{bmatrix}
      2 & 1 & 0 \\
     1 & 0 & -1 \\
      0 & -1 & -2 
    \end{bmatrix}$} & 
    
    \subcaptionbox{135$^{\circ}$}{$\begin{bmatrix}
      0 & 1 & 2 \\
      -1 &  0 &  1 \\
     -2 &  -1 &  0 
    \end{bmatrix}$} & 
    
    \subcaptionbox{225$^{\circ}$}{$\begin{bmatrix}
     -2 & -1 &  0 \\
      -1 &  0 &  1 \\
      0 &  1 &  2 
    \end{bmatrix}$} & 
    
    \subcaptionbox{315$^{\circ}$}{$\begin{bmatrix}
      0 & -1 & -2 \\
      1 &  0 & -1 \\
      2 &  1 &  0 
    \end{bmatrix}$} \\
    
  \end{tabular}
  }
  \caption{\textbf{Standard directional convolutional masks suitable for min-plus method.}}
  \label{t1}
\end{table}
This modification ensures that only the most dominant intensity variations contribute to edge detection, filtering out minor fluctuations and preserving prominent transitions. Unlike classical convolution, which considers all weighted contributions, tropical convolution emphasizes only the most significant local changes, thus effectively detecting sharper edges while suppressing noise.

\begin{example}
    Comparison Between Classical and Tropical Convolution
Consider a \(3 \times 3\) image patch $I =
\begin{bmatrix}
3 & 5 & 2 \\
6 & 1 & 4 \\
7 & 2 & 3
\end{bmatrix}$ and a convolution kernel $K =
\begin{bmatrix}
0 & 1 & 0 \\
1 & -4 & 1 \\
0 & 1 & 0
\end{bmatrix}.$

Applying classical convolution

\begin{equation}
O(x,y) = \sum_{i,j} K(i,j) \cdot I(x-i, y-j),
\end{equation}

which, for the center pixel \( I(1,1) = 1 \), results in

\begin{align*}
O(1,1) &= (0 \cdot 3) + (1 \cdot 5) + (0 \cdot 2) + (1 \cdot 6) + (-4 \cdot 1) + (1 \cdot 4) + (0 \cdot 7) + (1 \cdot 2) + (0 \cdot 3) \\
 &= 5 + 6 - 4 + 4 + 2\\
 &= 13.
\end{align*}

For tropical convolution, the operation is computed as

\begin{equation}
O_T(x,y) = \min_{i,j} \left[ K(i,j) + I(x-i, y-j) \right].
\end{equation}

\[
O_T(1,1) = \min(3,6,2,7,-3,5,7,3,3).
\]

\[
O_T(1,1) = -3.
\]
\end{example}

This demonstrates that tropical convolution selects the most prominent intensity drop, making it more effective at detecting edges while ignoring less relevant pixel variations.

Edges in an image correspond to sharp intensity changes, which can be detected by computing pixel-wise differences. Traditional gradient-based edge detection approximates these differences using finite difference operators, typically in the form

\begin{equation}
\delta I(x,y) = \left| I(x+1,y) - I(x,y) \right| + \left| I(x,y+1) - I(x,y) \right|.
\end{equation}

Instead of summing absolute differences, tropical edge detection redefines the gradient operator as:

\begin{equation}
\delta_T I(x,y) = \min \left( I(x+1,y) - I(x,y), I(x,y+1) - I(x,y) \right),
\end{equation}

where the minimum function captures the strongest downward intensity change. This formulation ensures that only the most dominant intensity drop is preserved, preventing weaker variations from influencing the edge detection process.

\begin{example}
    Consider the \(3 \times 3\) intensity matrix: $I =
\begin{bmatrix}
8 & 7 & 5 \\
6 & 4 & 3 \\
5 & 3 & 2
\end{bmatrix}.$

First, compute the horizontal differences:

\[
I(x+1,y) - I(x,y) =
\begin{bmatrix}
(6 - 8) & (4 - 7) & (3 - 5) \\
(5 - 6) & (3 - 4) & (2 - 3)
\end{bmatrix}
=
\begin{bmatrix}
-2 & -3 & -2 \\
-1 & -1 & -1
\end{bmatrix}.
\]

Next, compute the vertical differences

\[
I(x,y+1) - I(x,y) =
\begin{bmatrix}
(7 - 8) & (5 - 7)\\
(4 - 6) & (3-4)\\
(3 - 5) & (2 - 3)
\end{bmatrix}
=
\begin{bmatrix}
-1 & -2 \\
-2 & -1 \\
-2 & -1
\end{bmatrix}.
\]

Due to the difference in dimensions between the horizontal and vertical difference matrices, we either apply padding to match the original $3\times3$ size or restrict computation to the overlapping $2\times2$ region. In this case, we apply zero-padding to maintain consistent dimensions for further processing. Applying the tropical min operator

\[
\delta_T I(x,y) =
\min\left(
\begin{bmatrix}
-2 & -3 & -2 \\
-1 & -1 & -1 \\
 0 & 0 & 0
\end{bmatrix}
,
\begin{bmatrix}
-1 & -2 & 0 \\
-2 & -1 & 0\\
-2 & -1 & 0
\end{bmatrix}\right)
=
\begin{bmatrix}
-2 & -3 & -2 \\
-2 & -1 & -1 \\
-2 & -1 & 0
\end{bmatrix}.
\]
\end{example}
The negative values in this final matrix indicate strong edges, with larger negative values corresponding to more prominent edges.

We compute the intensity difference in the horizontal direction $I(x+1,y) - I(x,y)$ and vertical direction $I(x,y+1) - I(x,y)$ then we take the minimum of these two differences, ensuring we capture the strongest downward intensity change.

Let's analyze how it behaves
\begin{itemize}
    \item If \( I(x+1,y) \gg I(x,y) \), it means that the next pixel is significantly brighter than the current one and the difference \( I(x+1,y) - I(x,y) \) is positive and large. Since the tropical operator selects the minimum value, this large positive difference is ignored. As a result, no edge is detected because there is no strong downward transition in intensity.
    \item If \( I(x+1,y) \ll I(x,y) \), it indicates that the next pixel is much darker than the current one and the difference \( I(x+1,y) - I(x,y) \) is negative, representing a sharp intensity drop. Since the tropical operator takes the minimum value, this negative difference is selected. Consequently, an edge is detected at this location.
\end{itemize}

To enhance edge detection, preprocessing is performed using edge-preserving bilateral filtering followed by anisotropic shock filtering. Bilateral filtering smoothens the image while retaining sharp boundaries and is mathematically represented as  

\begin{equation}  
I'(x,y) = \frac{1}{W_p} \sum_{x_i \in S} I(x_i) f_r(||I(x_i) - I(x)||) g_s(||x_i - x||),  
\end{equation}  

where \( f_r \) and \( g_s \) define range and spatial weight functions, respectively. This process removes noise while preserving edge details. Further enhancement is performed using shock filtering, where the intensity gradient \( \nabla I(x,y) \) sharpens edge features according to  

\begin{equation}  
I_s(x,y) = I(x,y) - \lambda \cdot \nabla I(x,y),  
\end{equation}  

where \( \lambda \) is a weight controlling sharpening strength.  

To capture edges at multiple scales, image resizing is applied to generate lower-resolution (\( I_{\text{small}} \)) and higher-resolution (\( I_{\text{large}} \)) representations. This multi-scale approach ensures robust edge detection across varying spatial frequencies and is mathematically expressed as  

\begin{equation}  
I_{\text{small}} = \text{resize}(I, s), \quad I_{\text{large}} = \text{resize}(I, l),  
\end{equation}  

where \( s < 1 \) reduces resolution and \( l > 1 \) enhances detail preservation.  

Edges are then extracted using \textit{tropical max-plus convolution}, which replaces traditional Sobel and Prewitt filters with multi-directional convolutional kernels, ensuring robust edge detection in all orientations. The edge response is computed using  

\begin{equation}  
O_T^{\max}(x,y) = \max_{i,j} \left[ K(i,j) + I(x-i, y-j) \right],  
\end{equation}  

where the max operation ensures that only the strongest edge responses are retained. The filtering is performed using four directional masks, for ease we call it as 4 kernels and a six directional mask as 6 kernels, where as eight directional masks in Table[\ref{t1}] termed as 8 kernels for conventional use.

\[
h_1 =
\begin{bmatrix}
1 & 2 & 1 \\
0 & -4 & 0 \\
-1 & -2 & -1
\end{bmatrix}, \quad
h_2 =
\begin{bmatrix}
1 & 0 & -1 \\
2 & -4 & -2 \\
1 & 0 & -1
\end{bmatrix},
\]

\[
h_3 =
\begin{bmatrix}
-1 & -2 & -1 \\
0 & -4 & 0 \\
1 & 2 & 1
\end{bmatrix}, \quad
h_4 =
\begin{bmatrix}
-1 & 0 & 1 \\
-2 & -4 & 2 \\
-1 & 0 & 1
\end{bmatrix}.
\]

These filters extract edges in horizontal, vertical, and diagonal orientations, producing edge maps  

\begin{equation}  
E_T = \max \left( O_T^{\max}(h_1), O_T^{\max}(h_2), O_T^{\max}(h_3), O_T^{\max}(h_4) \right).  
\end{equation}  

To further refine detected edges, Hessian filtering is employed, leveraging second-order derivatives to enhance structural details. The Hessian matrix is defined as  

\begin{equation}  
H =
\begin{bmatrix}
\frac{\partial^2 I}{\partial x^2} & \frac{\partial^2 I}{\partial x \partial y} \\
\frac{\partial^2 I}{\partial y \partial x} & \frac{\partial^2 I}{\partial y^2}
\end{bmatrix},  
\end{equation}  

where the eigenvalues of \( H \) indicate curvature strength, enabling the detection of highly structured edges.  

To ensure robustness, wavelet shrinkage is applied to enhance contrast and suppress noise artifacts. The wavelet response is computed as  

\begin{equation}  
W_T(x,y) = | \nabla E_T(x,y) | \cdot (1 - E_T(x,y)),  
\end{equation}  

where wavelet transformation enhances edge response while reducing interference from non-edge regions.  

Edges are then binarized using adaptive thresholding, ensuring locally optimal segmentation:  

\begin{equation}  
T(x,y) = \frac{1}{N} \sum_{x',y' \in N} E_T(x',y'),  
\end{equation}  

where \( N \) is the local neighborhood. The final edge map is given by  

\begin{equation}  
E_{\text{final}}(x,y) =  
\begin{cases}  
1, & \text{if } E_T(x,y) \geq T(x,y), \\  
0, & \text{otherwise}.  
\end{cases}  
\end{equation}  

To ensure continuous and smooth edges, morphological thinning is applied via  

\begin{equation}  
E_{\text{thinned}} = \text{thin}(E_{\text{final}}),  
\end{equation}  

and area-based filtering removes small artifacts while preserving meaningful structures.  

Finally, image quality metrics are computed to assess edge performance. The contrast ratio is defined as  

\begin{equation}  
C = \frac{\sigma_{\text{edges}}}{\sigma_{\text{background}}},  
\end{equation}  

where \( \sigma \) denotes standard deviation. Correlation between detected and reference edges is given by  

\begin{equation}  
R = \frac{\sum (I_{\text{edges}} - \mu_1)(I_{\text{groundtruth}} - \mu_2)}{\sigma_1 \sigma_2},  
\end{equation}  

while energy and entropy provide further insights into image structure:  

\begin{equation}  
E = \sum P(x,y)^2,  
\end{equation}  

\begin{equation}  
H = -\sum P(x,y) \log P(x,y).  
\end{equation}  

By incorporating tropical convolution, Hessian filtering, wavelet shrinkage, and adaptive thresholding, the proposed method improves edge sharpness, reduces sensitivity to noise, and maintains computational efficiency. These characteristics make it a promising alternative to traditional edge detection approaches, particularly in low-contrast and structurally complex images.

\section{Experimental results}
For the evaluation of the proposed tropical geometry-based edge detection framework, we utilized a \(3 \times 3\) pixel window, which is computationally efficient and theoretically optimal. Unlike larger window sizes, such as \(5 \times 5\), which may introduce excessive blur and over-exposure due to the involvement of more metric parameters, the \(3 \times 3\) window ensures a balance between detail preservation and edge enhancement. The built-in images available in MATLAB, resized to \(300 \times 400\), were used to conduct our experiments. The key objective of our methodology is to enhance edge detection quality while maintaining structural continuity and computational efficiency. Our experimental setup was implemented using MATLAB (MathWorks) on a Windows 11 system with an 11th Gen Intel(R) Core(TM) i3-1115G4 processor running at 3.00 GHz. The computational efficiency of our method ensures that each test image is processed within a few seconds, demonstrating its real-time applicability in image processing tasks.

\begin{figure}
  \centering
  \begin{subfigure}[b]{0.161\linewidth}
    \centering
    Original\\    \includegraphics[width=\linewidth]{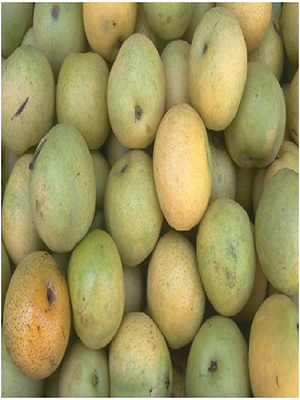}
    \caption*{EME =  0.202}
  \end{subfigure}
  \begin{subfigure}[b]{0.161\linewidth}
    \centering
    $TG^{\text{adapt tresh}}_{min-plus}$\\
    \includegraphics[width=\linewidth]{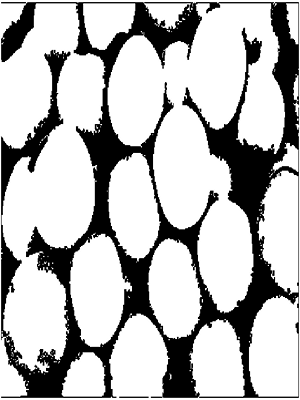}
    \caption*{EME = 6.491}
  \end{subfigure}
    \begin{subfigure}[b]{0.161\linewidth}
    \centering
    $TG^{\text{random tresh}}_{max-plus}$\\
    \includegraphics[width=\linewidth]{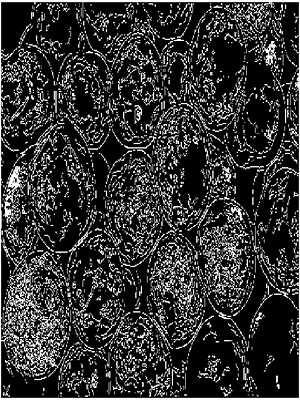}
    \caption*{EME =  8.383}
  \end{subfigure}
  \begin{subfigure}[b]{0.161\linewidth}
    \centering
    $TG^{\text{8 kernels}}_{min-plus}$\\
    \includegraphics[width=\linewidth]{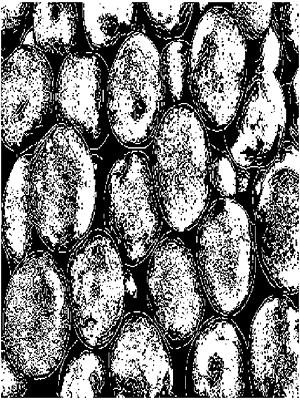}
    \caption*{EME = 7.829}
  \end{subfigure}
  \begin{subfigure}[b]{0.161\linewidth}
    \centering
    $TG^{\text{6 kernels}}_{max-plus}$\\
    \includegraphics[width=\linewidth]{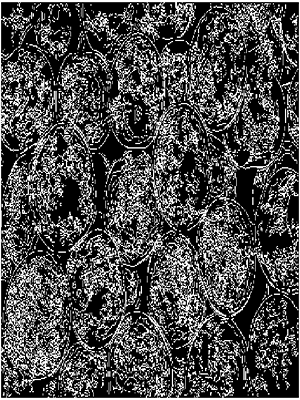}
    \caption*{EME = 8.397}
  \end{subfigure}
  \begin{subfigure}[b]{0.161\linewidth}
    \centering
    $TG^{\text{4 kernels}}_{max-plus}$\\
    \includegraphics[width=\linewidth]{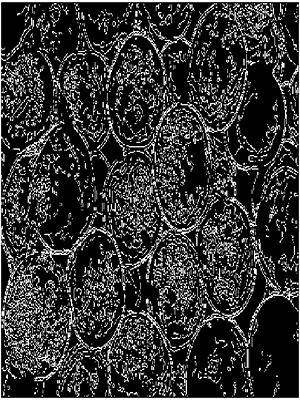}
    \caption*{EME =  8.380}
  \end{subfigure}
  \caption{Edge detection results using min-plus and max-plus algebraic operations}
  \label{f1}
\end{figure}

\begin{figure}
  \centering
  \begin{subfigure}[b]{0.2\linewidth}
    \includegraphics[width=\linewidth]{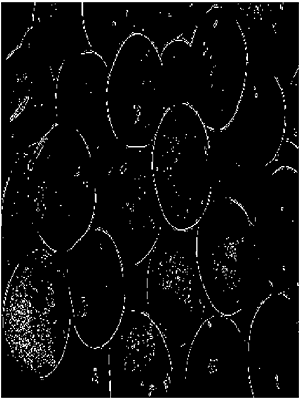}
    \caption*{Robert}
  \end{subfigure}
    \begin{subfigure}[b]{0.2\linewidth}
    \includegraphics[width=\linewidth]{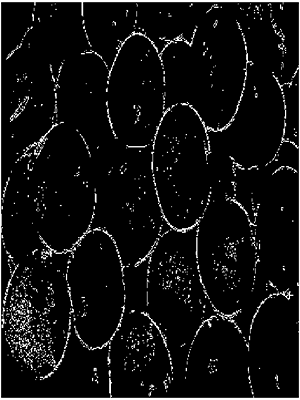}
    \caption*{$\text{TG}_{\text{Robert}}$}
  \end{subfigure}
  \begin{subfigure}[b]{0.2\linewidth}
    \includegraphics[width=\linewidth]{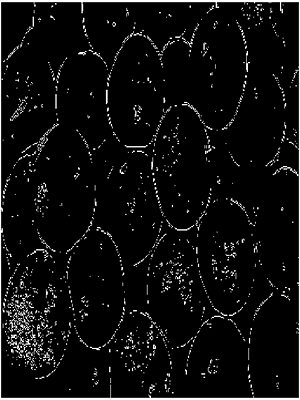}
    \caption*{Prewitt}
  \end{subfigure}
    \begin{subfigure}[b]{0.2\linewidth}
    \includegraphics[width=\linewidth]{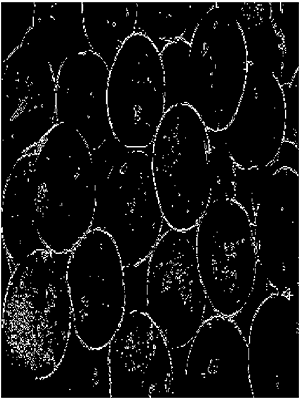}
    \caption*{${\text{TG}_{\text{Prewitt}}}$}
  \end{subfigure}
      \begin{subfigure}[b]{0.2\linewidth}
    \includegraphics[width=\linewidth]{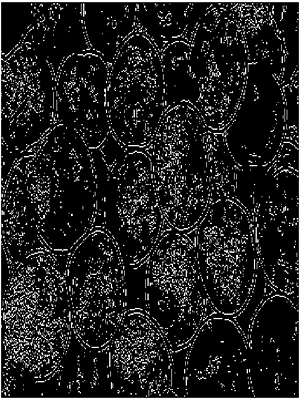}
    \caption*{Log}
  \end{subfigure}
      \begin{subfigure}[b]{0.2\linewidth}
    \includegraphics[width=\linewidth]{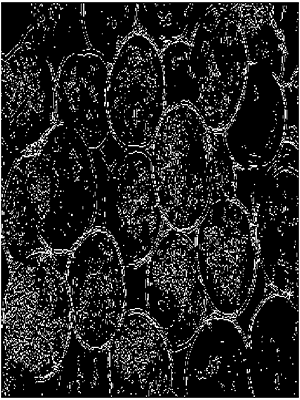}
    \caption*{${\text{TG}_{\text{Log}}}$}
  \end{subfigure}
      \begin{subfigure}[b]{0.2\linewidth}
    \includegraphics[width=\linewidth]{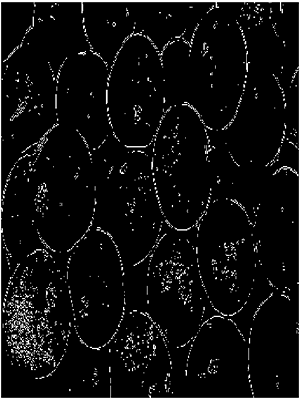}
    \caption*{Sobel}
  \end{subfigure}
        \begin{subfigure}[b]{0.2\linewidth}
    \includegraphics[width=\linewidth]{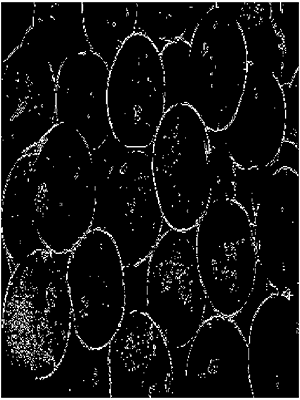}
    \caption*{${\text{TG}_{\text{Sobel}}}$}
  \end{subfigure}
  \caption{Visual comparison of edge detection using classical operators and TG-integrated operators}
  \label{various}
\end{figure}

\begin{figure}
  \centering
  \begin{subfigure}[b]{0.3\linewidth}
    \includegraphics[width=\linewidth]{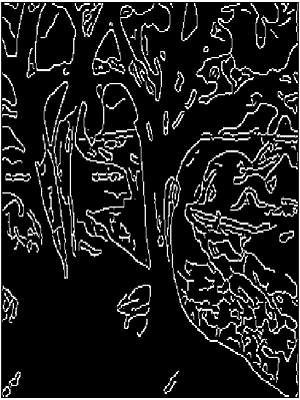}
    \caption*{Threshold = 0}
  \end{subfigure}
    \begin{subfigure}[b]{0.3\linewidth}
    \includegraphics[width=\linewidth]{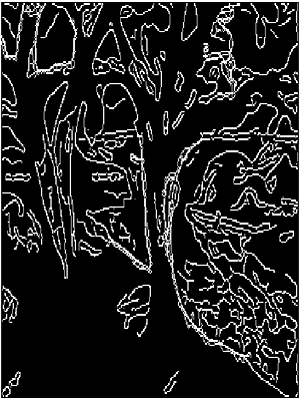}
    \caption*{Threshold = 0.5}
  \end{subfigure}
    \begin{subfigure}[b]{0.3\linewidth}
    \includegraphics[width=\linewidth]{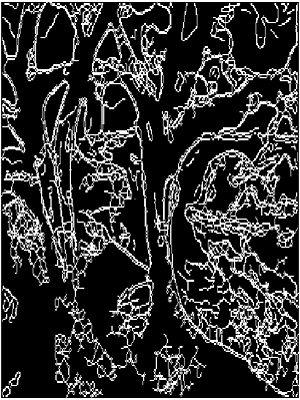}
    \caption*{Threshold = 1}
  \end{subfigure}
  \caption{Edge detection results on MATLAB's tree image with varying thresholding levels}
  \label{Threshold}
\end{figure}

\begin{figure}
\centering
\rule{\textwidth}{0.4pt} \\
\vspace{0.15cm}

\makebox[0.24\textwidth][c]{\textbf{MATLAB Images}}%
\hspace{0.005\textwidth}%
\makebox[0.24\textwidth][c]{\textbf{$TG^{\text{adapt tresh}}_{min-plus}$}}%
\hspace{0.005\textwidth}%
\makebox[0.24\textwidth][c]{\textbf{$TG^{\text{8 kernels}}_{min-plus}$}}%
\hspace{0.005\textwidth}%
\makebox[0.24\textwidth][c]{\textbf{$TG^{\text{4 kernels}}_{max-plus}$}}%
\\
\vspace{0.2cm} 
\rule{\textwidth}{0.4pt} \\

\begin{subfigure}[b]{0.24\textwidth}
\centering
\includegraphics[width=\linewidth]{pearsog.png}
\end{subfigure}%
\hspace{0.005\textwidth}%
\begin{subfigure}[b]{0.24\textwidth}
\centering
\includegraphics[width=\linewidth]{ptg5.png}
\end{subfigure}%
\hspace{0.005\textwidth}%
\begin{subfigure}[b]{0.24\textwidth}
\centering
\includegraphics[width=\linewidth]{ptg4.png}
\end{subfigure}%
\hspace{0.005\textwidth}%
\begin{subfigure}[b]{0.24\textwidth}
\centering
\includegraphics[width=\linewidth]{ptg1.png}
\end{subfigure}

\vspace{0.1cm} 
\begin{subfigure}[b]{0.24\textwidth}
\centering
\includegraphics[width=\linewidth]{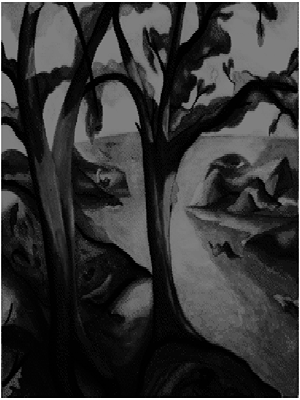}
\end{subfigure}%
\hspace{0.005\textwidth}%
\begin{subfigure}[b]{0.24\textwidth}
\centering
\includegraphics[width=\linewidth]{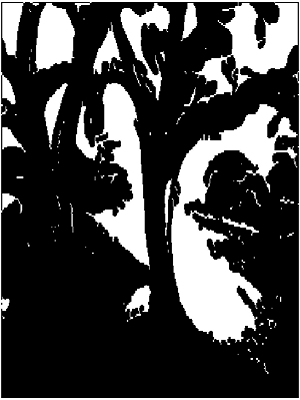}
\end{subfigure}%
\hspace{0.005\textwidth}%
\begin{subfigure}[b]{0.24\textwidth}
\centering
\includegraphics[width=\linewidth]{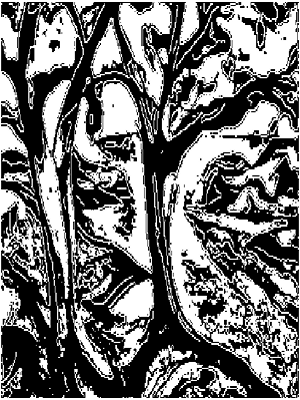}
\end{subfigure}%
\hspace{0.005\textwidth}%
\begin{subfigure}[b]{0.24\textwidth}
\centering
\includegraphics[width=\linewidth]{thres0.png}
\end{subfigure}

\vspace{0.1cm} 
\begin{subfigure}[b]{0.24\textwidth}
\centering
\includegraphics[width=\linewidth]{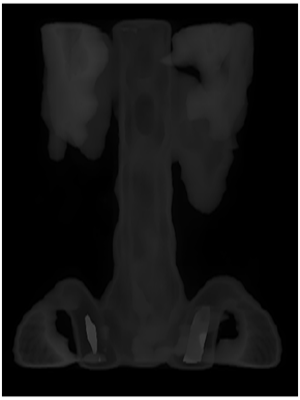}
\end{subfigure}%
\hspace{0.005\textwidth}%
\begin{subfigure}[b]{0.24\textwidth}
\centering
\includegraphics[width=\linewidth]{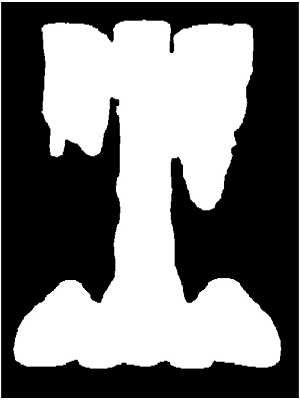}
\end{subfigure}%
\hspace{0.005\textwidth}%
\begin{subfigure}[b]{0.24\textwidth}
\centering
\includegraphics[width=\linewidth]{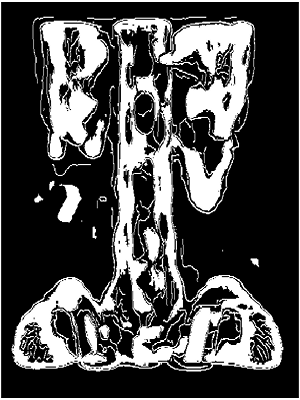}
\end{subfigure}%
\hspace{0.005\textwidth}%
\begin{subfigure}[b]{0.24\textwidth}
\centering
\includegraphics[width=\linewidth]{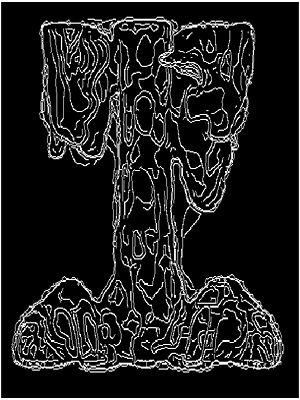}
\end{subfigure}

\vspace{0.1cm} 
\begin{subfigure}[b]{0.24\textwidth}
\centering
\includegraphics[width=\linewidth]{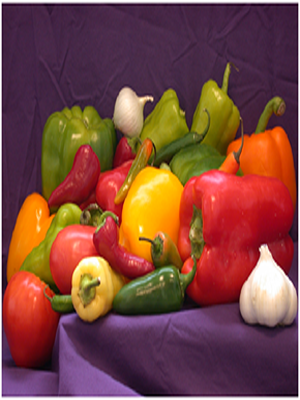}
\end{subfigure}%
\hspace{0.005\textwidth}%
\begin{subfigure}[b]{0.24\textwidth}
\centering
\includegraphics[width=\linewidth]{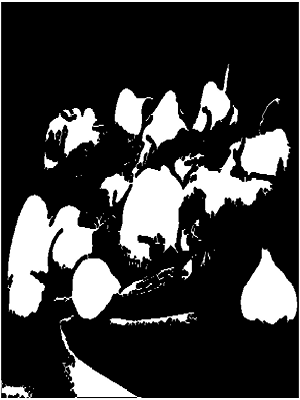}
\end{subfigure}%
\hspace{0.005\textwidth}%
\begin{subfigure}[b]{0.24\textwidth}
\centering
\includegraphics[width=\linewidth]{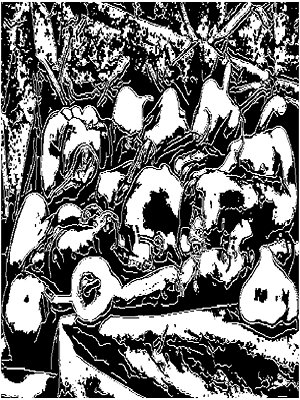}
\end{subfigure}%
\hspace{0.005\textwidth}%
\begin{subfigure}[b]{0.24\textwidth}
\centering
\includegraphics[width=\linewidth]{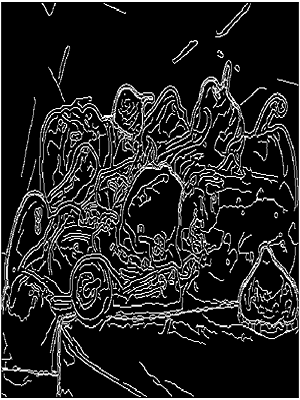}
\end{subfigure}

\vspace{0.1cm} 
\begin{subfigure}[b]{0.24\textwidth}
\centering
\includegraphics[width=\linewidth]{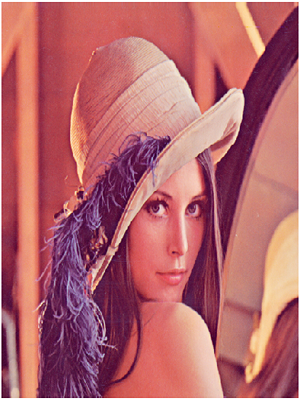}
\end{subfigure}%
\hspace{0.005\textwidth}%
\begin{subfigure}[b]{0.24\textwidth}
\centering
\includegraphics[width=\linewidth]{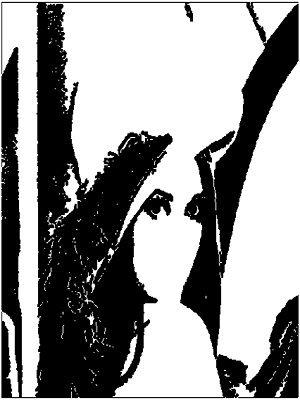}
\end{subfigure}%
\hspace{0.005\textwidth}%
\begin{subfigure}[b]{0.24\textwidth}
\centering
\includegraphics[width=\linewidth]{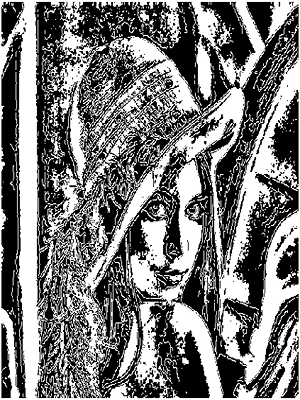}
\end{subfigure}%
\hspace{0.005\textwidth}%
\begin{subfigure}[b]{0.24\textwidth}
\centering
\includegraphics[width=\linewidth]{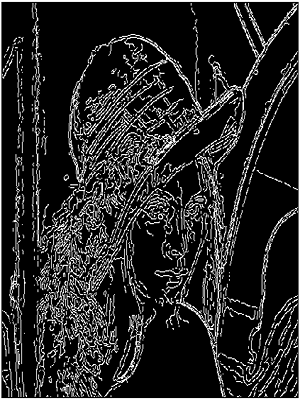}
\end{subfigure}

\caption{Visual comparison of edge detection methods applied to the MATLAB Images using $TG^{\text{adapt tresh}}_{min-plus}$, $TG^{\text{8 kernels}}_{min-plus}$ and $TG^{\text{4 kernels}}_{max-plus}$ methods.\label{vc}}
\end{figure}

\begin{table}
\centering
\renewcommand{\arraystretch}{1.5}
\setlength{\tabcolsep}{3pt}

\begin{adjustbox}{width=\textwidth}
\begin{tabular}{C{2.8cm} C{2cm} C{1.5cm} C{1.5cm} C{1.5cm} C{1.5cm} C{1.5cm} C{1.5cm} C{1.5cm} C{1.5cm} C{1.5cm} C{1.5cm}}
\hline
\textbf{Images} & \textbf{Methods} & \multicolumn{2}{c}{\textbf{Contrast}} & \multicolumn{2}{c}{\textbf{Correlation}} & \multicolumn{2}{c}{\textbf{Energy}} & \multicolumn{2}{c}{\textbf{Entropy}} & \multicolumn{2}{c}{\textbf{Homogeneity}} \\
\cmidrule(lr){3-12}
& & \textbf{Original} & \textbf{Edge} & \textbf{Original} & \textbf{Edge} & \textbf{Original} & \textbf{Edge} & \textbf{Original} & \textbf{Edge} & \textbf{Original} & \textbf{Edge} \\
\hline

\multirow{3}{*}{\includegraphics[width=0.45\linewidth]{pearsog.png}} 
& $TG^{\text{adapt tresh}}_{min-plus}$ &  & 1889.798 &  & 0.921 &  &  0.689 &  & 1.777 &  & 0.880 \\
& $TG^{\text{8 kernels}}_{min-plus}$ & 193.206 & 8558.781 & 0.945 & 0.642 & 0.021 & 0.265 & 7.264 & 4.928 & 0.198 & 0.397 \\
& $TG^{\text{4 kernels}}_{max-plus}$ &  & 7569.692 &  &  0.412 &  & 0.441 &  & 4.253 &  & 0.467 \\
\hline
\multirow{3}{*}{\includegraphics[width=0.45\linewidth]{treeog.png}} 
& $TG^{\text{adapt tresh}}_{min-plus}$ &  & 1673.927 &  &  0.927 &  & 0.724 &  &  1.544 &  & 0.901 \\
& $TG^{\text{8 kernels}}_{min-plus}$ & 371.594 & 6397.710 & 0.928 & 0.761 & 0.061 & 0.425 & 6.522 & 3.655 & 0.334 &  0.607 \\
& $TG^{\text{4 kernels}}_{max-plus}$ &  & 4207.212 &  & 0.591 &  & 0.733 &  &  2.177 &  & 0.764 \\
\hline
\multirow{3}{*}{\includegraphics[width=0.45\linewidth]{spineog.png}} 
& $TG^{\text{adapt tresh}}_{min-plus}$ &  & 751.458 &  &  0.976 &  & 0.677 &  & 1.294 &  & 0.959 \\
& $TG^{\text{8 kernels}}_{min-plus}$ & 365.181 & 3495.766 & 0.890 & 0.852 & 0.367 & 0.605 & 4.257 & 2.571 & 0.783 & 0.769 \\
& $TG^{\text{4 kernels}}_{max-plus}$ &  & 4439.425 &  & 0.525 &  & 0.689 &  & 2.549 &  & 0.712 \\
\hline
\multirow{3}{*}{\includegraphics[width=0.45\linewidth]{peppersog.png}} 
& $TG^{\text{adapt tresh}}_{min-plus}$ &  & 1422.335 &  &  0.936 &  & 0.741 &  & 1.428 &  & 0.912 \\
& $TG^{\text{8 kernels}}_{min-plus}$ & 267.145 & 7060.705 & 0.949 & 0.724 & 0.039 & 0.362 & 7.106 & 4.196 & 0.361 & 0.529 \\
& $TG^{\text{4 kernels}}_{max-plus}$ &  & 4970.418 &  & 0.525 &  & 0.664 &  & 2.753 &  & 0.688 \\
\hline
\multirow{3}{*}{\includegraphics[width=0.45\linewidth]{lenaog.png}} 
& $TG^{\text{adapt tresh}}_{min-plus}$ &  & 1689.341 &  &  0.944 &  &  0.635 &  & 1.802 &  & 0.895 \\
& $TG^{\text{8 kernels}}_{min-plus}$ & 370.480 & 8564.131 & 0.926 & 0.653 & 0.017 & 0.312 & 7.511 & 4.548 & 0.199 & 0.459 \\
& $TG^{\text{4 kernels}}_{max-plus}$ &  & 6388.476 &  & 0.405 &  & 0.594 &  & 3.229 &  & 0.610 \\

\hline
\end{tabular}
\end{adjustbox}

\caption{\textbf{Computational values of various quality metrics (contrast, correlation, energy, entropy, and homogeneity) between Original Images, $TG^{\text{adapt tresh}}_{min-plus}$, $TG^{\text{8 kernels}}_{min-plus}$ and $TG^{\text{4 kernels}}_{max-plus}$}}
\label{metrics}
\end{table}
The proposed framework adapted to five classical operators—Roberts, Prewitt, Sobel, LoG, and Canny—as illustrated in Figure~[\ref{various}]. Among the evaluated variants, the Tropical Canny (\(\text{TG}_\text{C}\)) method demonstrated improved sensitivity to inner boundaries, while the Tropical LoG (\(\text{TG}_\text{LoG}\)) yielded consistent results in low-noise, low-texture scenarios.  Figure~[\ref{f1}] illustrates the Edge Mean Entropy (EME)\cite{Agaian2000} for all well performing variant of TG. Additionally, the $TG^{\text{4 kernels}}_{max-plus}$ variant exhibited flexible behavior across threshold values, as shown in Figure~[\ref{Threshold}], producing continuous and well-connected edge structures in varying image conditions.

Quantitative assessments using standard quality metrics—contrast, entropy, homogeneity, correlation, and energy—are presented in Table~[\ref{metrics}]. These results suggest that the tropical framework is effective in preserving structural information and enhancing edge contrast across both grayscale and color images. Figure~[\ref{vc}] provides a visual comparison of edge detection results on MATLAB sample images using three tropical variants: $TG^{\text{adapt tresh}}{min-plus}$, $TG^{\text{8 kernels}}{min-plus}$, and $TG^{\text{4 kernels}}_{max-plus}$. The results illustrate how tropical algebra can be integrated with conventional edge detection pipelines to enhance their performance. The adaptive threshold-based variant yields clear boundary separation in low-noise scenarios. The multi-kernel min-plus variant captures fine-grained edge textures and internal object structures with high detail retention. The max-plus variant emphasizes structural smoothness and continuity, effectively connecting edge fragments. Together, these methods demonstrate the adaptability of the tropical framework to different edge detection tasks, complementing and extending the capabilities of classical operators such as Canny and LoG.

Overall, the results indicate that the tropical geometry-based framework offers a viable alternative to conventional edge detection techniques, particularly in scenarios involving low contrast or complex structures. The algebraic operations used in the method contribute to consistent edge representation with reduced sensitivity to noise, supporting its potential applicability in a range of image analysis tasks.

\section{Conclusion and Future Work}
This study presented a tropical geometry-based edge detection framework that reformulates convolution and gradient operations using min-plus and max-plus algebra. The method enhances edge sharpness, structural continuity, and robustness to noise by emphasizing dominant intensity transitions through algebraic transformations. It has been evaluated using grayscale and color images, showing effective integration with classical operators such as Canny, Sobel, Prewitt, and LoG.

Among the tested variants, the Tropical Canny (\(\text{TG}_\text{C}\)) and Tropical LoG (\(\text{TG}_\text{LoG}\)) demonstrated improved detection of inner boundaries and consistent performance in low-noise, low-texture environments. In addition, three core tropical formulations—$TG^{\text{adapt tresh}}_{\text{min-plus}}$, $TG^{\text{8 kernels}}_{\text{min-plus}}$, and $TG^{\text{4 kernels}}_{\text{max-plus}}$—illustrate the adaptability of the framework. Each offers distinct benefits, such as enhanced boundary definition, high detail retention, or improved structural smoothness, depending on the image characteristics.

Quantitative assessments using metrics such as contrast, entropy, and Edge Mean Entropy (EME) support the framework’s capability to produce well-defined and continuous edge maps with minimal computational overhead. While the results are promising, the current implementation relies on fixed kernel sizes and static parameter settings, which may limit adaptability across different image domains. The method has been tested on standard images, and broader validation on diverse and large-scale datasets remains an area for future work. Additionally, although the tropical operations are computationally efficient in MATLAB, further optimization may be needed for real-time or embedded applications.

Overall, the framework offers a scalable and mathematically grounded alternative to traditional edge detection approaches. Future extensions will explore adaptive parameter tuning, deployment in resource-constrained environments, and integration with deep learning architectures to support advanced feature extraction and classification tasks.

\subsection*{Data Availability}
The dataset utilized in this study are publicly accessible, ensuring transparency and reproducibility of our results.

\section*{Declarations}

\subsection*{Ethics approval and consent to participate}
Not applicable

\subsection*{Consent for publication}
Not applicable

\subsection*{Competing interests}
The authors declare no competing interests.

\end{document}